\documentclass{elsarticle}
\usepackage{epsfig}
\usepackage{float}
\usepackage[ansinew]{inputenc}
\usepackage{amsmath,amsfonts,latexsym,amssymb,amsbsy}
\usepackage{graphicx}
\usepackage{amsmath,amssymb}
\usepackage{pgfplots}
	\usetikzlibrary{plotmarks}
		\usetikzlibrary{external}
		\pgfplotsset{compat=1.3}
		\newlength\figurewidth 
		\newlength\figureheight
		\tikzset{external/force remake=true} 
\usepackage{array}
\usepackage{booktabs}
\usepackage{xcolor}
\definecolor{marin}{rgb}   {0.,   0.3,   0.7}
\definecolor{rouge}{rgb}   {0.8,   0.,   0.}
\definecolor{sepia}{rgb}   {0.8,   0.5,   0.}
\usepackage[colorlinks,citecolor=marin,linkcolor=rouge,
            bookmarksopen,
            bookmarksnumbered
           ]{hyperref}
\usepackage{epsf}

\newcommand{\R}{\mathbb{R}}

\newcommand{\F}{\mathcal{F}}

\newcommand{\e}{\ensuremath{\mathrm{e}}}

\newcommand{\cO}{{\mathcal O}}

\numberwithin{equation}{section}

\newcommand{\QED}{\mbox{}\hfill \raisebox{-0.2pt}{\rule{5.6pt}{6pt}\rule{0pt}{0pt}}
          \medskip\par}

\begin{document}
\begin{frontmatter}
\title{Numerical solution of the Burgers' equation with high order splitting methods }

\author[rvt,els]{M. Seydao\u{g}lu\corref{cor1}}
\ead{muasey@imm.upv.es}
\author[focal]{U. Erdo\u{g}an}
\ead{utku.erdogan@usak.edu.tr}
\author[els]{T. \"{O}zi\c{s}}
\ead{turgut.ozis@ege.edu.tr}

\cortext[cor1]{Corresponding author}
\address[rvt]{Department of Mathematics, Faculty of Art and Science,
Mu\c{s} Alparslan University, 49100 Mu\c{s}, Turkey.}
\address[focal]{Department of Mathematics, Faculty of Art and Science, U\c{s}ak University, 64200, U\c{s}ak, Turkey.}
\address[els]{Department of Mathematics, Faculty of Science, Ege University, 35100, Izmir, Turkey}

\begin{abstract}
In this work, high order splitting methods have been used for calculating the numerical solutions of the Burgers' equation in one space dimension with periodic and Dirichlet boundary conditions. However, splitting methods with real coefficients of order higher than two necessarily have negative coefficients and can not be used for time-irreversible systems, such as Burgers equations, due to the time-irreversibility of the Laplacian operator. Therefore, the splitting methods with complex coefficients and extrapolation methods with real and positive coefficients have been employed. If we consider the system as the perturbation of an exactly solvable problem(or can be easily approximated numerically), it is possible to employ highly efficient methods to approximate Burgers' equation. The numerical results show that the methods with complex time steps having one set of coefficients real and positive, say $a_i\in\mathbb{R}^+$ and $b_i\in\mathbb{C}^+$, and high order extrapolation methods derived from a lower order splitting method produce very accurate solutions of the Burgers' equation.

\end{abstract}


\end{frontmatter}

\section{Introduction} 
In this work we consider the numerical integration of the one-dimensional Burgers' equation
\begin{equation}\label{eq:Burger}
  \frac{\partial u}{\partial t}=\frac{\partial }{\partial  x} (\frac{u^{2}}{2})+\nu \frac{\partial^{2} u}{\partial x^{2}},  \qquad u(x,0)=u_{0}(x).
\end{equation}
where $t$ and $x$ represent time and distance in x-direction respectively,  and $\nu > 0$ is the kinematic viscosity
parameter related to the Reynolds number $R=1/\nu$, was first given its steady state solutions by Bateman \cite{bateman15}. It was noticed later, by Burgers \cite{burger48}, this equation is used in a mathematical modeling of turbulence. It is exploited in a different applied fields, such as in the context of gas dynamics, in a modeling shock waves, traffic flows and continuous stochastic process. On the other hand, Burgers' equation can be solved exactly for several initial data by Hopf-Cole transformation and the solutions can be represented as a Fourier series expansion \cite{hopf50,Cole51}. For different initial data, exact solutions are listed in \cite{ben72zman}. Since the exact solutions are available for some limiting cases, i.e. some set of initial functions, it  is  natural to consider Burgers' equation for testing the performance of a numerical solver. Thus, many interesting papers on the numerical solution of Burgers' equation  based on the finite difference, the finite element, the boundary element and the spectral methods exist in the literature \cite{cal82smith,eva84abd,mitta93sign,ozis96ozdes,kutbah99esen,kut04esen,has05salhos,w08lia,oz09utku,tu05eskut}. 

Jain and Raja \cite{jain79raja} have seperated Burgers' equation in two sub models, namely, the convection and the diffusion part, then used the finite difference method to solve each model problem successively and called as "splitting-up technique". Similar to this strategy, Jain and Holla \cite{jain78holla} have used the cubic spline method and called "two-time-level splitting". Furthermore, the new algorithms along with the cubic spline method is proposed in \cite{jain92shansing}, which treats Burgers' equation as three sub problems, referred as "three-time-level splitting". We refer to \cite{saka07dag} for time and space splitting techniques applied to the Burgers' equation and the modified Burgers' equation along quintic B-spline collocation procedure.

However, Godunov and Strang splitting methods have been applied to the Generalized Korteweg-De Vries Equations(KdV) and convergence properties of the both methods have been analyzed in \cite{hold99karrisb}. In addition, a new analytical approach to the  Godunov and Strang splitting methods presented in \cite{hold11krtao} for the KdV equation and in \cite{hold13lubris} for the particular examples of the PDEs with Burgers' nonlinearity, by using requirement of a well-posedness theory in Sobolev spaces. 
A wide variety of techniques have been considered up to recently, which have been limited to low orders of accuracy in time due to the combination of nonlinearity and stiffness in Burgers' equation. To obtain high order accurate solutions of stiff PDEs, the fourth-order modified exponential time-differencing method is presented in \cite{kassam05treft}.

In this work we are interested in numerical integration of the Burgers' equation using high order splitting methods with complex coefficients and extrapolation methods. Such procedures could allow us to integrate irreversible systems. Splitting methods with complex coefficients have been recently presented for the numerical integration of the autonomous and non-autonomous parabolic equations in \cite{blanes13oho,castella09smw,hansen09hos,blanes13mua}. Since the Burgers' equation involves small viscosity parameter $\nu$, we also consider splitting methods which designed for the perturbation of an exactly solvable problems (or can be easily approximated numerically). In usual, the presence of small $\nu$ is considered to be a numerical challenge. However utilizing perturbed system can take the advantage of smaller parameter.

It has been shown that splitting methods with real coefficients of order higher than two involve negative time steps \cite{blanes05otn,goldman96nos,sheng89slp,suzuki90fdo} and therefore this is undesirable for Burgers' equation. For this reason, we consider splitting methods with complex coefficients having positive real part and with real positive coefficients obtained by applying extrapolation techniques to a lower order splitting method. We use spectral methods for periodic boundary conditions, finite difference and Weighted Essentially Nonoscillatory (WENO) schemes for dirichlet boundary conditions as spatial discretization techniques. This allows us to adopt methods for system of ordinary differential equations (ODEs) to solve the partial differential equation (\ref{eq:Burger}) numerically.
\section{Splitting methods} 
Let us consider the nonlinear parabolic PDE of the form
\begin{eqnarray} \label{eq:para}
\frac{d u}{dt}=  \hat A(u(t))  + \hat B(u(t)), \qquad u(0)=u_0,
\end{eqnarray}
$u(x,t)\in\R^D$, and where the (possibly unbounded) operators
$\hat A$, $\hat B$ and $\hat A+\hat B$ generate $C^0$ semi-groups
for positive $t$ over a finite or infinite Banach space. 

For the sake of simplicity, we write the non-linear equation in the
(apparently) linear form
\begin{eqnarray} \label{eq:paralie}
\frac{d u}{dt}=  Au(t)  + Bu(t),
\end{eqnarray}
where $A,B$ are the Lie operators associated to $ \hat A, \hat B$,
i.e.
\begin{eqnarray} \label{eq:lieoprt}
  A \equiv \hat A(u(t)) \frac{\partial}{\partial u}, \qquad
  B \equiv \hat B(u(t)) \frac{\partial}{\partial u}
\end{eqnarray}
which act on functions of $u(t)$. The formal solution of the (\ref{eq:paralie}) is given by $u(t)=\e^{t(A+B)}u_0$, which is a
short way to write
\[
  u(t)=\e^{t(A+B)}u_0 = \left.
  \sum_{k=0}^{\infty} \frac{t^k}{k!}\left(\hat A(u(t)) \frac{\partial}{\partial u} +
  \hat B(u(t)) \frac{\partial}{\partial u} \right)^k u \
  \right|_{u=u_0}.
\]
The framework of the splitting method for solving numerically (\ref{eq:paralie}) is to decompose the solutions into the exactly (or numerically) solvable  two subproblems
\begin{eqnarray} \label{eq.Split_aut}
 \frac{d u}{dt}=  A u(t) \qquad \mbox{ and } \qquad
 \frac{d u}{dt}=  B u(t).
\end{eqnarray}
and treat them seperately.
It is possible to replace notation of the formal solution of the (\ref{eq:paralie}), $\e^{t(A+B)}$, by the flow map $\Phi_{t}^{( A + B)}$ in the nonlinear case. Let denote by $\Phi^{A}_{h}, \Phi^{B}_{h}$ ( or $\e^{hA}$, $\e^{hB}$) the exact $h$-flows for each problem in
(\ref{eq.Split_aut}) (and for a sufficiently small time step, $h$). Then 
the simplest method within this class is 
\begin{equation}
L_{h} = \Phi^{A}_{h} \circ \Phi^{B}_{h} \qquad \mbox{ or } \qquad L_{h} = \Phi^{B}_{h} \circ  \Phi^{A}_{h},
\end{equation}
which is known as the first order {\em Lie-Trotter
splitting} method. A splitting method $\Psi_{h}$ has order $p$, if $\Psi_{h} = \Phi_{h}^{( A + B)}+ \cO(h^{p+1})$. If one composes Lie-Trotter method and its adjoint $L^{\ast}_{h} = \Phi^{B}_{h} \circ \Phi^{A}_{h}$ with a half time step, one obtains second order time symmetric methods
\begin{equation} \label{eq:strang2}
 S_{h}= \Phi^{A}_{h/2} \circ \Phi^{B}_{h} \circ \Phi^{A}_{h/2}, 
\end{equation}
\begin{equation} \label{eq:strang3}
 S_{h}= \Phi^{B}_{h/2} \circ \Phi^{A}_{h} \circ \Phi^{B}_{h/2},
\end{equation}
which are referred as {\em Strang splitting} based on the pattern $"ABA"$ and $"BAB"$ respectively. For the nonlinear case high-order approximations  based on the pattern $"BAB"$
can be obtained
\begin{eqnarray} \label{eq:splittingmethod}
\Psi_{h} =  \Phi_{h b_1 }^{B} \circ \Phi_{h a_1 }^{A} \circ \cdots  \circ  \Phi_{h b_m}^{B} \circ \Phi_{h a_m }^{A} \circ \Phi_{h b_{m+1}}^{B}  ,
\end{eqnarray}
or, equivalently
\begin{eqnarray} \label{eq:splittingmethod1}
\Psi_{h} = \e^{h b_{m+1} B} \, \e^{h a_m A} \,\cdots \, \e^{h b_2
B} \, \e^{h a_1 A} \, \e^{h b_1 B},
\end{eqnarray}
where m is usually referred as number of steps of the method, and $a_{1}, ..., a_{m}$ and $b_{1}, ..., b_{m}$ are real or complex numbers(to be chosen) depending on the desired order of method \cite{creutz89hhm,suzuki90fdo,yoshida90coh}. If one takes $b_1=b_{m+1}=0$, then one transforms the pattern $"BAB"$ into the pattern $"ABA"$
with a different computational cost. Notice that the difference between the order of the indices in equations (\ref{eq:splittingmethod}) and (\ref{eq:splittingmethod1}) is appeared by the Lie-derivative action and called as \textit{Vertauschungssatz} \cite{hairer06gni}.

Splitting methods with real coefficients of order
greater than two necessarily contain some negative coefficients. In general, the semi-groups are not well defined for negative time steps. Positivity requirement on the coefficients prevents the use of splitting methods of order greater than two when operators $A$ and $B$ generates a $C^0$ semi-group of propagators. Independently, the results in \cite{castella09smw} and \cite{hansen09hos} resolve the open question for the existence of the splittings methods of order two in the context of semi-groups. Additionally, the authors of \cite{castella09smw} and \cite{hansen09hos} derived a new broad class of splitting methods with complex coefficients of order 3 to 14 by a method which is referred "triple-jump composition procedure", and presented theoretical error bounds in the linear case . At least formally, all results given in \cite{castella09smw,hansen09hos} can instantly be stretched out to nonlinear case by replacing all exponentials with the corresponding nonlinear flows. The numerical results for the autonomous and non-autonomous nonlinear cases can be found in \cite{castella09smw,blanes13mua}.

However, as already mentioned, one can consider using splitting methods designed for near integrable system to solve Burgers' equation because of the term involving small viscosity parameter $\nu$. Therefore, $\hat A(u(t))$ is dominant part, i.e $\|\hat B(u(t))\|\ll\|\hat A(u(t))\|$. Furthermore, one can write equations (\ref{eq:paralie}) as
\begin{eqnarray} \label{eq:pertsmall}
\frac{d u}{dt}=  (A + \epsilon B)u(t),
\end{eqnarray}
where $\epsilon$ is small parameter, i.e $|\varepsilon|\ll 1$.

Let us consider the symmetric second-order methods (\ref{eq:strang2}) with $B$ replaced by $\epsilon B$ for the autonomous case in the exponential form
\begin{eqnarray} \label{eq:strlocerr}
\e^{\frac{h}{2} A} \, \e^{\epsilon h B} \,\e^{\frac{h}{2} A}=\e^{h(A+\epsilon B)-\frac{1}{24} h^{3}[A,[A,\epsilon B]]+\frac{1}{12} h^{3}[\epsilon B,[A,\epsilon B]]+\,\cdots \,}.
\end{eqnarray}
The Lie bracket (or commutator) $[.,.]$, which is defined as $[A,B]=AB-BA$ for A and B denoted in (\ref{eq:lieoprt}), arises from the application of the Baker-Campbell-Hausdorff formula\cite{hairer06gni} to the left-hand side of (\ref{eq:strlocerr}). It is easy to see the local error of the above method is of order $\cO(\epsilon h^{3}+\epsilon^{2} h^{3})$ and originates in the commutators at third order in h,
i.e $\epsilon [A,[A,B]]$ and $\epsilon^{2}[ B,[A,B]]$. Essentially, a small parameter $\epsilon$ is considered to be $\epsilon \ll\ h$ or at least $\epsilon \approx h$. Thus, one can cancel the dominant error terms in $\epsilon$ rather than in $h$ for the general composition (\ref{eq:splittingmethod1}) with $B$ replaced by $\epsilon B$ and built methods which take advantage of this relevant property by choosing the coefficients $a_{i},b_{i}$.
An $m$-stage symmetric $BAB$ compositions given by
\begin{eqnarray} \label{eq:BAB}
\Psi(h) = \e^{h b_{m+1}\varepsilon B} \, \e^{h a_m A} \,\cdots \,
\e^{h b_2\varepsilon B} \, \e^{h a_1 A} \, \e^{h b_1\varepsilon
B},
\end{eqnarray}
with $a_{m+1-i}=a_i, \ b_{m+2-1}=b_i, \ i=1,2,\ldots$,  and $ABA$
compositions are given by
\begin{eqnarray} \label{eq:ABA}
\Psi(h) = \e^{h a_{m+1} A} \, \e^{h b_m \varepsilon B} \,\cdots \,
\e^{h a_2 A} \, \e^{h b_1\varepsilon B} \, \e^{h a_1 A},
\end{eqnarray}
with $a_{m+2-i}=a_i, \ b_{m+1-1}=b_i, \ i=1,2,\ldots$. 
 
In these cases, the dominant error terms can be read as $\cO(\epsilon h^{p_{1}})$ where $p_{1}$ is considered as relatively large values. Then, one can take into account the small parameter $\epsilon$ in the accuracy of the desired splitting methods. Let $(p_{1},p_{2})$ be an effective order of a method  with $p_{1}\geq p_{2}$ that yields the local error $\cO(\epsilon h^{p_{1}+1}+\epsilon^{2} h^{p_{2}+1})$. Some methods of order $(p_{1},2)$ for $p_{1}\leq 10$ with all coefficients $a_{i}, b_{i}$ are positive and methods of order $(p_{1},4)$ for $p_{1}=6,8$ presented in \cite{mclah95lan}. In \cite{blanes13nfo}, the order conditions of the symmetric splitting methods of a given generalized order $(p_{1},2)$ for $p_{1}\geq 4$ and  $(p_{1},4)$ for $p_{1}\geq 6$ are presented by using Lyndon multi-indices and some schemes of order $(p_{1},4)$ for $p_{1}=8,10$ are obtained by solving corresponding order conditions. 

On the other hand, due to the cost of a multiplication, splitting methods with complex coefficients are computationally more costly than with real coefficients. To reduce the computational cost some new high order methods for which only one set of coefficients complex valued are proposed in \cite{castella09smw,blanes13mua}. Additionally, the effective order of $(6,4)$ method is obtained in \cite{blanes13mua} with the coefficients $a_{i}$ being positive and real valued, whereas the coefficients $b_{i}$ being complex valued with positive real part. For non-autonomous perturbed systems, this last method is most efficient and stable.

\subsection{Splitting methods for Burgers' equation}
Considering initial value problem for the viscous Burgers' equation of the form
\begin{equation}\label{eq:vBurger}
  u_{t}=(\frac{u^{2}}{2})_{x}+\nu u_{xx},  \qquad u(x,0)=u_{0}(x).
\end{equation}
We now describe the framework of the splitting method, that is, to solve in succession the conservation law
\begin{equation}\label{eq:conservation}
  u_{t}=(\frac{u^{2}}{2})_{x},  \qquad u(x,0)=u_{0}(x),
\end{equation}
and the diffusion equation 
\begin{equation}\label{eq:diffusion}
  u_{t}=\nu u_{xx},  \qquad u(x,0)=u_{0}(x).
\end{equation}

Let us denote by $ \Phi^{A}_{h}, \Phi^{B}_{h} $ the maps associated to the exact
solution (or a sufficiently accurate numerical approximation) of the (\ref{eq:conservation}), (\ref{eq:diffusion}) respectively. Then, we approximate the solution of (\ref{eq:vBurger}) as
\begin{eqnarray} \label{eq:solnt}
u(x,h)=\Psi_{h}u_{0}(x). 
\end{eqnarray}
where $\Psi_{h}$ given by (\ref{eq:splittingmethod})(or \ref{eq:splittingmethod1}) for small time step $h$.

The readers can be suspicious about composing the solutions from (\ref{eq:conservation}) and (\ref{eq:diffusion}) since (\ref{eq:diffusion}) always produces smooth solutions while (\ref{eq:conservation}) results in discontinuous shock solution within a certain time interval depending on initial profile. We refer the readers two important research papers on this issue. Holden et.al \cite{hold99karrisb} showed that if the initial data are sufficiently regular, the Strang splitting method converges to the smooth solution of full equation provided that the splitting step size for Burgers' step (\ref{eq:conservation}) is kept under control. Another study by Holden et.al \cite{hold13lubris} proved that splitting solution converges to the weak solution of the full equation assuming that the splitting procedure is convergent. 
\subsection{Methods for the Diffusion equation}
In this section we consider the methods which have been used to approximate the diffusion equation (\ref{eq:diffusion}). 
\subsubsection{Fast Fourier transform }
Consider the equation (\ref{eq:diffusion}) for $x\in[0,2\pi]$, $t>0$, with periodic boundary conditions. Then, we can write the solution $u(x,t)$ as
\begin{eqnarray} \label{eq:fourier}
u(x,t)=\sum^{\infty}_{k=-\infty} \hat u_{k}(t)e^{ikx}, 
\end{eqnarray}
where $\hat u_{k}$ are Fourier coefficients of the initial function. The formula for the discrete Fourier transform (DFT) is 
\begin{eqnarray} \label{eq:dft}
\hat u_{k}(t)=h\sum^{N}_{j=1} \hat u_{j}(t)e^{-ikx_{j}}, \qquad  k=-N/2+1, \ldots , N/2,
\end{eqnarray}
where the spacing of the grid point is $h=\frac{2\pi}{N}$ and the inverse DFT is given by
\begin{eqnarray} \label{eq:idft}
u_{j}(t)=\frac{1}{2\pi}\sum^{N/2}_{k=-N/2+1} \hat u_{k}(t)e^{ikx_{j}}, \qquad  j=1, \ldots , N.
\end{eqnarray}
Above DFT can be computed by the Fast Fourier Transform (FFT) \cite{trefethen20}. Now DFT for the equation (\ref{eq:diffusion}) reads
\begin{eqnarray} \label{eq:dftburg}
\frac{d}{dt}\hat u_{k}(t)=-\nu k^{2} \hat u_{k}(t), \qquad \hat u_{k}(0)= \hat u_{0k},
\end{eqnarray}
of which solution can be written as follows
\begin{equation}\label{eq:slndftburg}
 \hat u_{k}(t)= e^{-\nu k^{2} t}\hat u_{0k}, \qquad  k=-N/2+1, \ldots , N/2. \nonumber
\end{equation}

\subsubsection{Finite difference method}
To apply splitting schemes in case of homogeneous  Dirichlet boundary conditions, we employ fourth order finite difference method. If one discretizes equation (\ref{eq:diffusion}) in space
$$x_j= j (\delta x), \qquad j=1,\ldots,D \quad \mbox{ with } \quad \delta x = 1/D,
$$
then one obtains the differential equation
\begin{equation} \label{eq:problem0}
\frac{dU}{dt} = BU,
\end{equation}
where $U=(U_1, \ldots,U_D)=(u_1, \ldots,u_D) \in \R^D$. The
Laplacian $\Delta$ has been approximated by the sparse matrix $B$ of size
$D\times D$ given by following discrete diffusion operator

\begin{equation}
B=\dfrac{1}{12 (\delta x)^2}\left( \begin{array}{cccccccc}
45 & -154 & 214&-156&61&0&\hdots&0 \\
16 & -30 & 16&-1& 0&\hdots & &0 \\
-1&16&-30&16&-1&0 &\hdots &0\\
\vdots\\

0&\hdots&0&61&-156&214&-154&45 \end{array} \right)
\end{equation}
and the solution of the (\ref{eq:problem0}) can be easily computed by using Matlab's \textit{expm}-routine.
\subsection{Methods for the Conservation Law}
In this section we consider the methods which have been used to approximate the Conservation Law (\ref{eq:conservation}). 

\subsubsection{Fast Fourier transform }
If one uses DFT as described for the diffusion equation then the equation (\ref{eq:conservation}) becomes
\begin{eqnarray} \label{eq:dftcons}
\frac{d}{dt}\hat u_{k}(t)=\frac{i}{2} k \hat u^{2}_{k}(t), \qquad \hat u_{k}(0)= \hat u_{0k}.
\end{eqnarray}
where $ k=-N/2+1, \ldots , N/2$. On the other hand we can write above equation as
\begin{eqnarray} \label{eq:dftcons1}
\frac{d}{dt}\hat u_{k}(t)=\frac{i}{2} k \F((\F^{-1}(\hat u_{k}(t)))^{2}), \qquad \hat u_{k}(0)= \hat u_{0k}.
\end{eqnarray}
$\F$ is the Fourier transform operator \cite{trefethen20}. For our numerical experiment, we use fourth order runge-kutta scheme for nonlinear ODE (\ref{eq:dftcons1}) which can be done in a similar way given by \cite[p 111]{trefethen20}.

\subsubsection{WENO finite differences }
  Weighted Essentially  Non-oscillatory  (WENO) schemes for spatial discretization are  proved to be very successful in the numerical treatment of convection dominated problems. The main idea of WENO finite difference  is to use an adaptive interpolation or reconstruction procedure based on the local smoothness of the numerical solution to get high order accuracy and oscillation free behavior near discontinues. For the conservation law 
\begin{equation} \label{conservation}
u_t+(f(u))_x=0 \nonumber
\end{equation} 
 The derivative $(f(u))_x$ is approximated  by
\begin{equation} \label{eq:flux}
(f(u))_x \vert _{x=x_j} \approx \frac{1}{\Delta x} (\hat{f}_{j+1/2}-\hat{f}_{j-1/2})
\end{equation}
where $\hat{f}_{j+1/2}$ is the numerical flux. In case of $f'(u)\geq 0$ ,the numerical flux of the fifth order WENO  finite difference  is given as follows \cite{FluidWeno}
\begin{equation}
\hat{f}_{j+1/2}=w_1\hat{f}_{j+1/2} ^{(1)}+w_2\hat{f}_{j+1/2} ^{(2)}+w_3\hat{f}_{j+1/2} ^{(3)} \nonumber
\end{equation}
Indeed the $\hat{f}_{j+1/2} ^{(i)}$ s are third order fluxes given by
\begin{eqnarray}
\hat{f}_{j+1/2} ^{(1)}&=&\frac{1}{3} f(u_{j-2})-\frac{7}{6} f(u_{j-1})+\frac{11}{6} f(u_{j}) \nonumber \\
\hat{f}_{j+1/2} ^{(2)}&=&\frac{-1}{6} f(u_{j-1})+\frac{5}{6} f(u_{j})+\frac{1}{3} f(u_{j+1}) \nonumber \\
\hat{f}_{j+1/2} ^{(3)}&=&\frac{1}{3} f(u_{j})+\frac{5}{6} f(u_{j+1})-\frac{1}{6} f(u_{j+2}) \nonumber
\end{eqnarray}
The  non linear weights in \ref{eq:flux} are given by
\begin{eqnarray}
w_i=\frac{\tilde{w}_i}{\sum _{k=1} ^{3} \tilde{w}_k} \nonumber \\
\tilde{w}_k=\frac{\gamma _k}{(\epsilon + \beta _k)^2} \nonumber
\end{eqnarray}
where the linear weights $\gamma_1=\frac{1}{10}$, $\gamma_2=\frac{3}{5}$ and $\gamma_3=\frac{3}{10}$. $\epsilon$ is taken $10^{-6}$  in actual computations. The smoothness indicators $\beta _k$ are listed below
\begin{eqnarray}
\beta_1&=& \frac{13}{12}\left( f(u_{j-2})-2 f(u_{j-1})+ f(u_{j}) \right)^2 +\frac{1}{4}\left( f(u_{j-2})-4 f(u_{j-1})+3 f(u_{j}) \right)^2 \nonumber \\
\beta_2&=& \frac{13}{12}\left( f(u_{j-1})-2 f(u_{j})+ f(u_{j+1}) \right)^2 +\frac{1}{4}\left( f(u_{j-1})- f(u_{j+1}) \right)^2  \nonumber\\
\beta_3&=& \frac{13}{12}\left( f(u_{j})-2 f(u_{j+1})+ f(u_{j+2}) \right)^2 +\frac{1}{4}\left(3 f(u_{j})-4 f(u_{j+1})+f(u_{j+2}) \right)^2  \nonumber
\end{eqnarray}
For detailed derivation of interpolation relations and reconstruction process we refer to the review article \cite{Weno}. In numerical experiments, the ODEs system arising from WENO discretization of conservation law is solved by explicit fourth order Runge-Kutta. 
\begin{table}[H]
\centering
\caption{Coefficients for several splitting schemes with the pattern "BAB".}
\label{tab.1}
 {
\begin{tabular}{ll}
\hline\hline
\multicolumn{2}{c}{The splitting method of effective order $(6,2)$: ML$(6,2)$} \\
\hline\hline
$ b_1= 1/12  $ & $a_1= (5-\sqrt{5})/10$\\
$ b_2= 5/12  $ & $a_2=1/\sqrt{5}$\\
$ b_3=b_2, \ b_4=b_1$ & $a_3=a_1$ 
\\
\hline\hline
\multicolumn{2}{c}{The 4-stage fourth-order method: RC$4$ } \\
\hline\hline
$ b_1= 1/10 - i/30$ & $a_1=a_2= a_3=a_4= 1/4 $\\
$ b_2= 4/15 + 2i/15$ \\
$ b_3=4/15-i/5$ \\
$ b_4=b_2, \ b_5=b_1$
\\
\hline\hline
\multicolumn{2}{c}{The optimized 4-stage fourth-order method: O$4$ } \\
\hline\hline
$ b_1= 0.060078275263542 + 0.060314841253379i$ & $a_1= 0.1859688195991091314 $\\
$ b_2= 0.270211839133611 - 0.152903932291162i$ & $a_2=  0.3140311804008908686 $\\
$ b_3=0.339419771205694 + 0.185178182075567i$ & $a_3=a_2, \ a_4=a_1$ \\
$ b_4=b_2, \ b_5=b_1$
\\
\hline\hline
\multicolumn{2}{c}{The optimized 4-stage fourth-order method: SM$4$ } \\
\hline\hline
$ b_1= 0.018329102861074364-0.10677008344599524i$ & $a_1= 0.13505265889288437$\\
$ b_2= 0.2784394345454581+0.20041452008768607i$ & $a_2= 0.36494734110711563$\\
$ b_3= 0.40646292518693505-0.18728887328338165i$ & $a_3=a_2, \ a_4=a_1$ \\
$ b_4=b_2, \ b_5=b_1$
\\ \hline\hline
\multicolumn{2}{c}{The splitting method of effective order $(6,4)$: SM(6,4) } \\
\hline\hline
$ b_1= 0.05753968253968254 - 0.007886748775536424i$ & $a_1= a_2=a_3=a_4=a_5=a_6=1/6$\\
$ b_2= 0.20476190476190473 + 0.04732049265321855i$ \\
$ b_3= 0.16309523809523818 - 0.11830123163304637i$  \\
$ b_4=0.14920634920634912 + 0.15773497551072851i$  \\
$ b_5=b_3, \ b_6=b_2, \ b_7=b_1$
\\ \hline\hline
\end{tabular}
}
\end{table}
\section{Numerical Results}
\begin{figure}[H]
\centering
	\pgfplotsset{every axis plot/.append style={line width=1.0pt, mark size=2pt},
		tick label style={font=\footnotesize},
		every axis/.append style={%
		minor x tick num=1,
		minor y tick num=4,
		minor z tick num = 2,
		scale only axis, 
		font=\footnotesize
		}
	}
	\setlength\figurewidth{.50\textwidth}
\setlength\figureheight{.45\textwidth}
	\includegraphics[width=.48\textwidth]{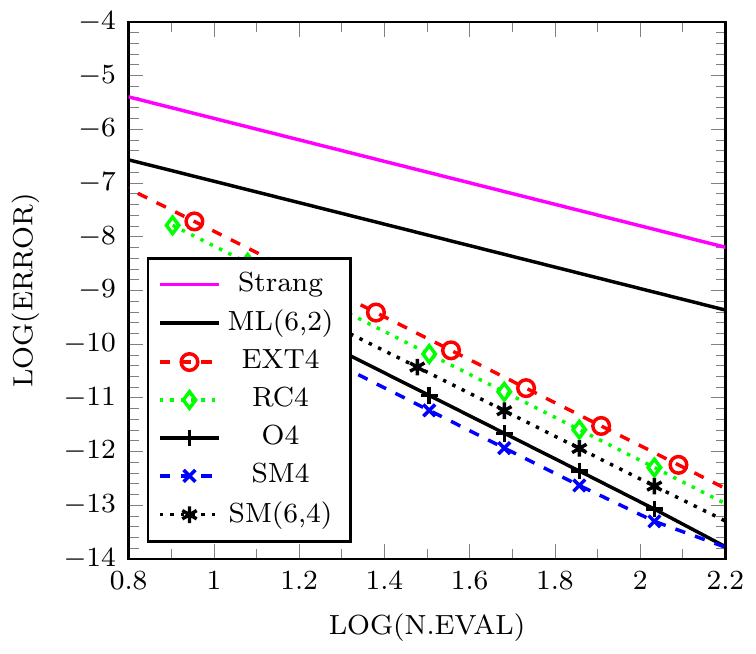}
	\includegraphics[width=.48\textwidth]{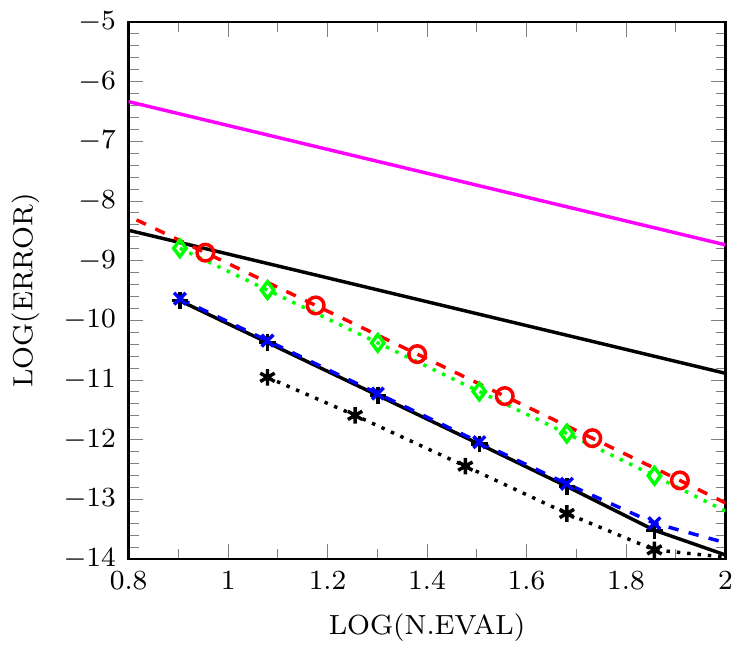}
\caption{\label{fig:ex1} Error versus number of evaluations of $\Phi^{A}_{h}$ for the numerical
    integration in Example 1 at $t=1$ for $\nu = 0.03$
     (left panel) and $\nu = 0.003$ (right panel).  }
\end{figure}

In this section, we numerically illustrate the performance of the different higher-order splitting methods, which are useful when highly accurate solutions of the one-dimensional problem (\ref{eq:Burger}) are sought. To overcome positivity requirements on the coefficients for the achieving second order barrier we first consider extrapolation methods
\begin{equation} \label{eq:extrpol}
 \Psi_{h}= \frac{4}{3} \Phi_{h/2} \circ \Phi_{h/2} -\frac{1}{3}\Phi_{h}.
\end{equation}
and
\begin{equation} \label{eq:extrpol6}
 \Psi_{h}= \frac{81}{40} \Phi_{h/3} \circ \Phi_{h/3} \circ \Phi_{h/3} -\frac{16}{15}\Phi_{h/2} \circ \Phi_{h/2}+\frac{1}{24}\Phi_{h}.
\end{equation}
If one takes the Strang splitting method  (\ref{eq:strang3}) as the basic method $\Phi_{h}$ with the exact flows, then one gets fourth-order method as  
\begin{equation} \label{eq:extrpol1}
 \Psi_{h}= \frac{4}{3} \Phi^{B}_{h/4} \circ \Phi^{A}_{h/2} \circ \Phi^{B}_{h/2} \circ \Phi^{A}_{h/2} \circ \Phi^{B}_{h/4} -\frac{1}{3}\Phi^{B}_{h/2} \circ \Phi^{A}_{h} \circ \Phi^{B}_{h/2},
\end{equation}
and sixth-order method as  
\begin{equation}
\begin{aligned} \label{eq:extrpol16}
 \Psi_{h} & = \frac{81}{40} \Phi^{B}_{h/6} \circ \Phi^{A}_{h/3} \circ \Phi^{B}_{h/3} \circ \Phi^{A}_{h/3} \circ \Phi^{B}_{h/3} \circ \Phi^{A}_{h/3} \circ \Phi^{B}_{h/6} \\
& -  \frac{16}{15} \Phi^{B}_{h/4} \circ \Phi^{A}_{h/2} \circ \Phi^{B}_{h/2} \circ \Phi^{A}_{h/2} \circ \Phi^{B}_{h/4} + \frac{1}{24}\Phi^{B}_{h/2} \circ \Phi^{A}_{h} \circ \Phi^{B}_{h/2},
\end{aligned}
\end{equation}
respectively.
We illustrate the results for the following schemes with real coefficients :

\begin{itemize}
\item {\bf Strang}: The second-order symmetric Strang splitting method (\ref{eq:strang3});
\item {\bf ML(6,2)}: The second-order symmetric splitting method built for perturbed systems in \cite{mclah95lan};
\item {\bf (EXT4)}: The fourth-order extrapolation method (\ref{eq:extrpol1});
\item {\bf (EXT6)}: The sixth-order extrapolation method (\ref{eq:extrpol16});
\end{itemize}
and we illustrate the results for the following schemes with complex coefficients and
$a_i\in\mathbb{R}^+$ :
\begin{itemize}
\item {\bf (RC4)}: The 4-stage fourth-order method from \cite{castella09smw};
\item {\bf (O4)}: The 4-stage fourth-order method built in \cite{blanes13oho},
 whose coefficients are available at
 \texttt{http://www.gicas.uji.es/Research/splitting-complex.html};
\item {\bf (SM4)}: The optimized 4-stage fourth-order method built for the perturbed systems in \cite{blanes13mua};
\item {\bf (SM(6,4))}: The 6-stage fourth-order method built for the perturbed systems in \cite{blanes13mua};
\end{itemize}
Coefficients of the above splitting schemes are given in Table~\ref{tab.1} for the convenience of the reader and have been considered to solve Burgers' equation with periodic boundary conditions and Dirichlet boundary conditions. The most appropriate methods are  symmetric BAB composition methods with all $a_{i}$ real and positive valued, $b_{i}$ complex valued having positive real part when solving Burgers' equation with periodic boundary conditions. If one considers spectral methods as a space discretization methods for periodic boundary conditions, this class of methods, namely $a_i\in\mathbb{R}^+$ and $b_i\in\mathbb{C}^+$ are stable and have less computational cost. For the Dirichlet case, this class of methods are not stable due to finite difference and WENO scheme, which have been used as spatial discretization techniques. In this case, we only use methods with real and positive time steps. On the other hand, for a given method $\Psi_{h}$ which is involve complex time steps, the numerical solutions $u_{n+1}$ computed by projection of the complex solutions to its real part after completing each time step, namely $u_{n+1}=\Re(\Psi_{h}u_{n})$.
\paragraph*{\bf Example 1} We consider the simulation of the Burgers' equation (\ref{eq:Burger}) with
\begin{equation} \label{eq:iintcond1}
 u(x,t=0)=\frac{1}{2}+\frac{1}{4}sin(x)
\end{equation}
and periodic boundary conditions in the space domain $[0, 2\pi]$. We take $\nu = 0.03$, $\nu = 0.003$ and the number of grid points as $N=512$ for Fourier spectral discretization in x. We compute the exact solution numerically by using fourth order Runge-Kutta methods based on the method of integrating factors given in \cite{trefethen20} for a sufficiently small time step. We measure the error of numerical solution at the end of the time integration in the infinity norm. In Fig.~\ref{fig:ex1}, we compare the accuracy of the splitting methods given in Table~\ref{tab.1} on the time interval $[0, 2\pi]$. We simulate the solution error versus the number of evaluations of $\Phi^{A}_{h}$ which usually requires the more costly computation for several step sizes. For all methods, we clearly observe the classical orders $p$ from the slopes of lines. Clearly, splitting methods with complex coefficients are slightly more accurate than lower order splitting methods with real coefficients and high order extrapolation methods. Furthermore, the standard methods are insensitive w.r.t. the small parameter $\nu$, whereas the splitting methods improve as $\nu$ decreases.

\begin{figure}[H]
\centering
	\pgfplotsset{every axis plot/.append style={line width=1.0pt, mark size=2pt},
		tick label style={font=\footnotesize},
		every axis/.append style={%
		minor x tick num=1,
		minor y tick num=4,
		minor z tick num = 2,
		scale only axis, 
		font=\footnotesize
		}
	}
	\setlength\figurewidth{.50\textwidth}
\setlength\figureheight{.45\textwidth}
	\includegraphics[width=.48\textwidth]{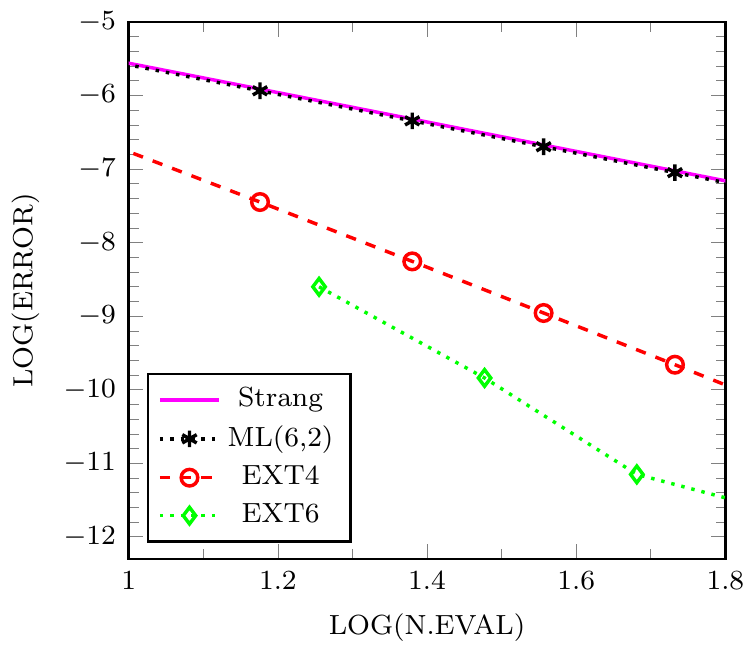}
	\includegraphics[width=.48\textwidth]{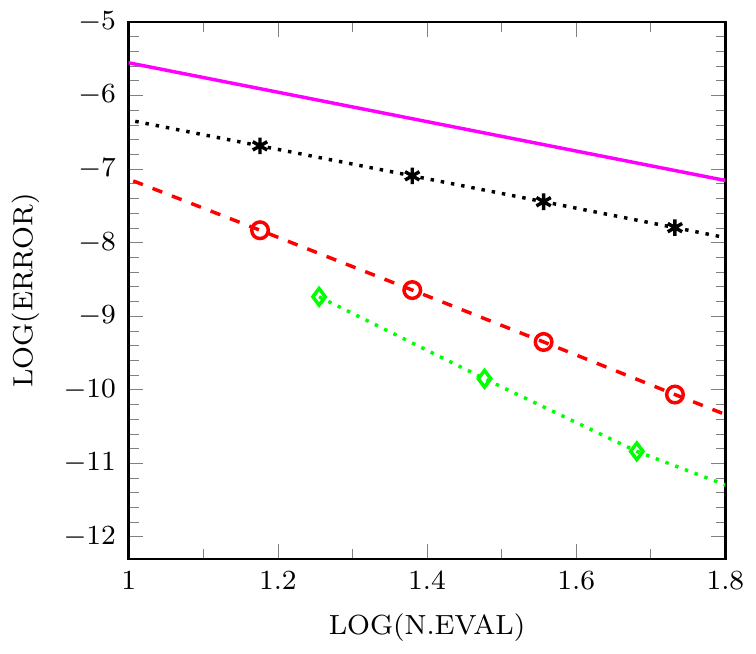}
		\includegraphics[width=.48\textwidth]{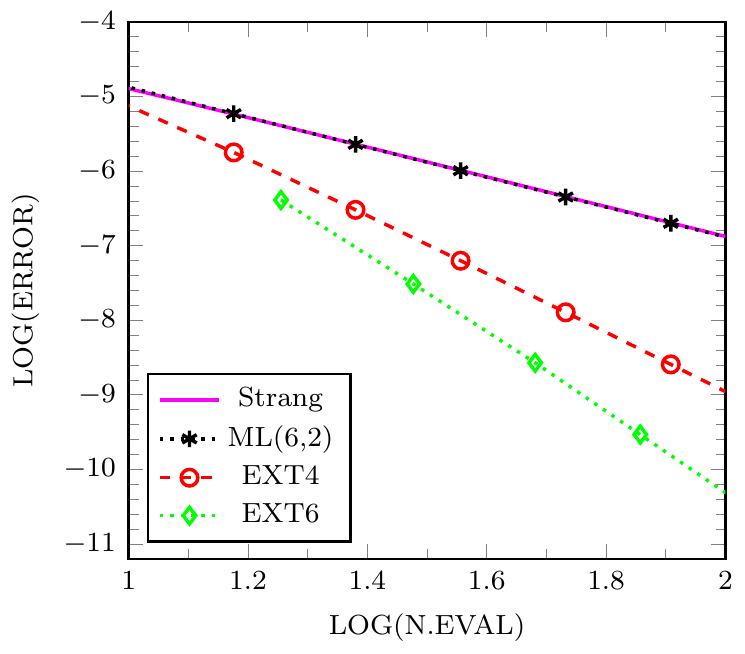}
		\includegraphics[width=.48\textwidth]{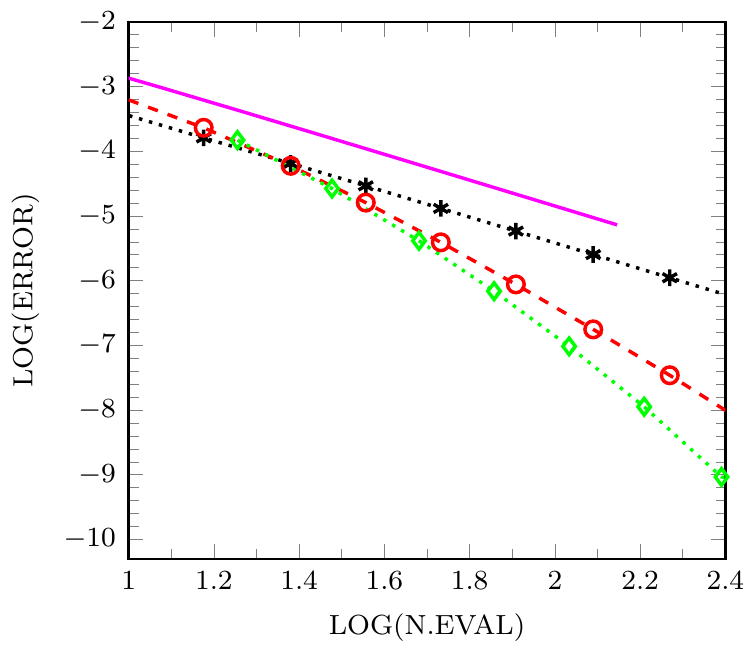}
\caption{\label{fig:ex2}  Error versus number of evaluations of $\Phi^{A}_{h}$ for the numerical
    integration in Example 2 at $t=1(Top),t=3(Bottom)$ for $\nu = 0.1$
 (left panel) and $\nu = 0.01$ (right panel).  }
\end{figure}

\paragraph*{\bf Example 2} The next test example corresponds the (\ref{eq:Burger}) on space domain $[0,1]$ with the initial condition
\begin{equation} \label{eq:iintcond2}
 u(x,t=0)= \frac{1}{5}sin(\pi x),
\end{equation}
and the following boundary conditions 
\begin{equation} \label{eq:boundcond1}
u(x=0,t)=u(x=1,t)=0, \qquad t > 0.  
\end{equation}
Using the Hopf-Cole transformation, the exact solution for this particular problem is given as follows 
\begin{equation} \label{eq:exactprob1}
 u(x,t)=2\nu \pi \frac{\sum^{\infty}_{n=1} c_{n} exp(-n^{2} \pi^{2} \nu t) n sin(n \pi x)}{c_{0}+ \sum^{\infty}_{n=1} c_{n} exp(-n^{2} \pi^{2} \nu t) cos(n \pi x)}, 
\end{equation}
where
\begin{eqnarray} \label{eq:exactprob11}
 c_{0} &=& \int^{1}_{0} exp\left\{-(10 \pi \nu)^{-1}\left[1-cos(\pi x)\right]\right\}dx, \nonumber \\
 c_{n} &=& 2 \int^{1}_{0} exp\left\{-(10 \pi \nu)^{-1}\left[1-cos(\pi x)\right]\right\} cos(n \pi x)dx \qquad (n=1,2,3 \ldots).  \nonumber
\end{eqnarray}
We take $\nu = 0.1$, $\nu = 0.01$ and the size of the discrete diffusion matrix $D=500$. We compute the infinity norm error of the numerical solution with respect to (\ref{eq:exactprob1}) at the final times $t=1$, $t=3$ by applying the compositions methods given in Table~\ref{tab.1}. The results can be seen in Fig.~\ref{fig:ex2}. As discussed in the paper \cite{hansen09hos}, error terms are in general not uniformly bounded on the interval $[0,T]$ for some positive $T$ in the infinite dimensional space when one imposes boundary conditions. Thus the convergence order is no longer guaranteed. For this reason, we observe severe order reductions in the experiments with Dirichlet boundary conditions. One clearly observes in Fig.~\ref{fig:ex2} that the extrapolation schemes are superior than lower order splitting schemes with real coefficients. However, even though the full orders are not obtained, high order splitting schemes produce considerable smaller errors than (6,2) and Strang splitting methods. The sixth order extrapolation scheme produce very accurate results among other schemes in this experiments. It is clear that the splitting method designed for perturbed system drastically improves when decreasing $\nu$.
\paragraph*{\bf Example 3} The last test example corresponds to the (\ref{eq:Burger}) on space domain $[0,1]$ with the initial condition
\begin{equation} \label{eq:iintcond3}
 u(x,t=0)= \frac{1}{2} x (1-x),
\end{equation}
and the same boundary conditions (\ref{eq:boundcond1}).
The exact solution for this particular problem is given by (\ref{eq:exactprob1}) with the coefficients 
\begin{eqnarray} \label{eq:exactprob21}
 c_{0} &=& \int^{1}_{0} exp\left\{-x^{2}(24 \nu)^{-1}(3-2x)\right\}dx, \nonumber \\
 c_{n} &=& 2 \int^{1}_{0} exp\left\{-x^{2}(24 \nu)^{-1}(3-2x)\right\} cos(n \pi x)dx \qquad (n=1,2,3 \ldots).  \nonumber
\end{eqnarray}
\begin{figure}[!t] 
\centering
	\pgfplotsset{every axis plot/.append style={line width=1.0pt, mark size=2pt},
		tick label style={font=\footnotesize},
		every axis/.append style={%
		minor x tick num=1,
		minor y tick num=4,
		minor z tick num = 2,
		scale only axis, 
		font=\footnotesize
		}
	}
	\setlength\figurewidth{.50\textwidth}
\setlength\figureheight{.45\textwidth}
	\includegraphics[width=.48\textwidth]{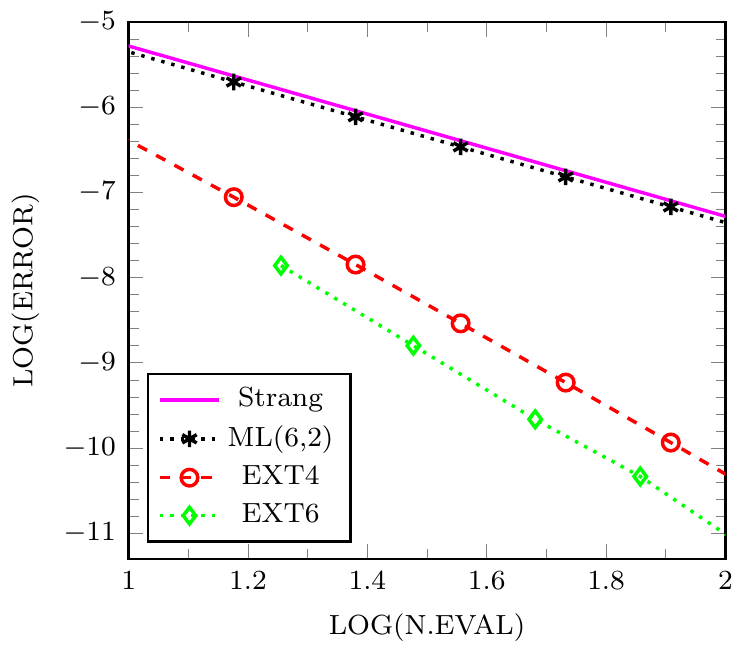}
	\includegraphics[width=.48\textwidth]{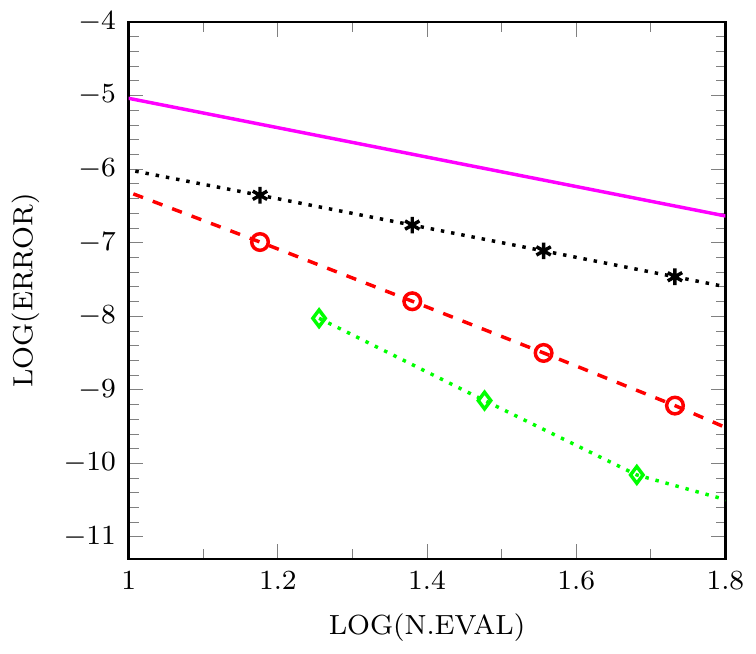}
	\includegraphics[width=.48\textwidth]{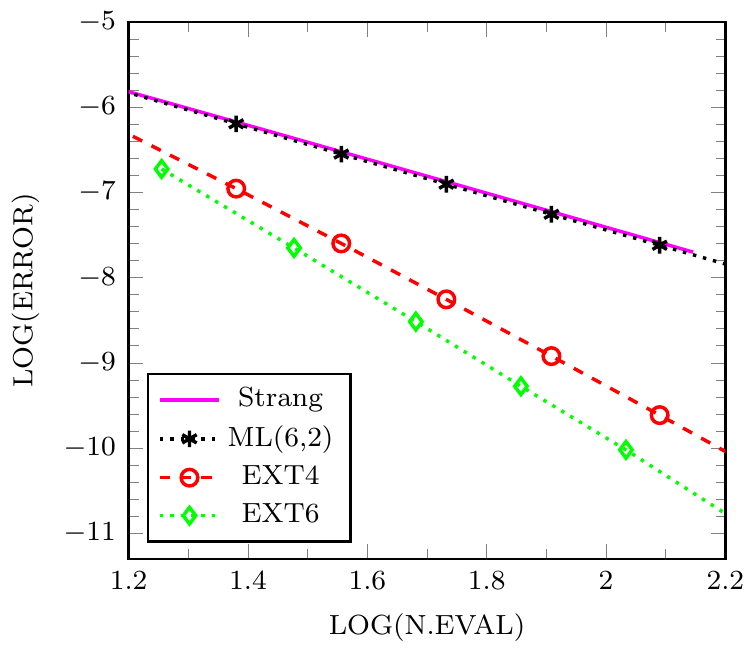}
	\includegraphics[width=.48\textwidth]{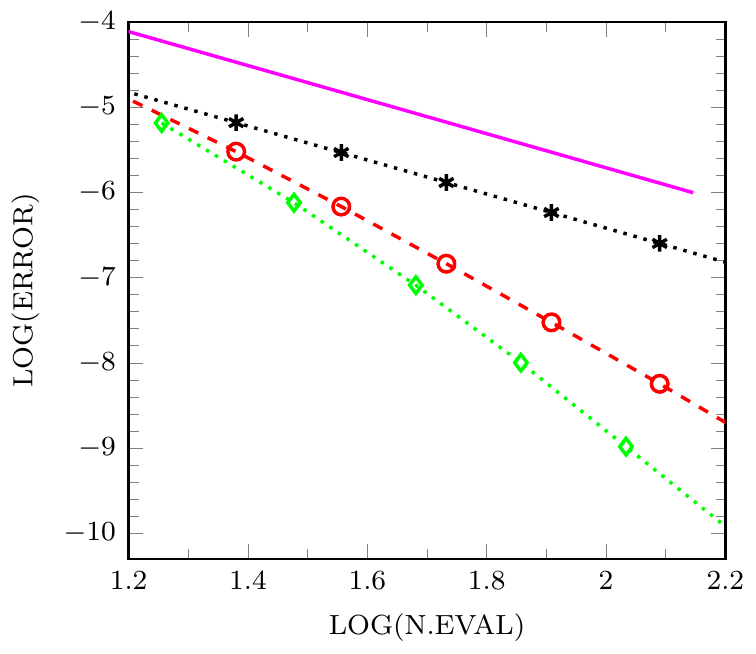}
\caption{\label{fig:ex3}  Error versus number of evaluations of $\Phi^{A}_{h}$ for the numerical
    integration in Example 2 at $t=1(Top),t=3(Bottom)$ for $\nu = 0.1$
 (left panel) and $\nu = 0.01$ (right panel).}
\end{figure}
We choose $\nu = 0.1$, $\nu = 0.01$, the size of the discrete diffusion matrix $D=500$ and compute the error at the final times $t=1$, $t=3$ by applying the same splitting schemes as in the Example 2. The results are collected in Fig.~\ref{fig:ex3} where the superiority of the extrapolation methods are manifest. Furthermore, the sensitivity w.r.t small parameter of the splitting methods which designed for perturbed systems also evident. 
\section{Conclusions}
We have considered the numerical integration of non-linear
Burgers' equations using high order splitting methods
with complex coefficients and real positive coefficients. Although there exists many high order accurate  numerical methods for pure diffusion and pure advection equations, designing a stable and efficient  method for singularly perturbed full PDEs is a challenging task. The suitable methods for diffusion and non linear advection could be applied subsequently without any changes through the higher order splitting methods. To overcome second order barrier of classical splitting algorithms with positive coefficients, many types of splitting procedures for Burgers' equation with periodic and Dirichlet boundary conditions are discussed through the paper. As alternatives to Lie Trotter and Strang schemes, higher order methods are derived by extrapolation and complex coefficients. Besides, effective spatial discretizations of the subequations are considered depending on the types of boundary conditions. In the numerical examples, the expected order reductions for the Burgers' equation with Dirichlet boundary conditions on bounded domains are reported. It is concluded that the methods designed for perturbed problems taking the advantage of small viscosity number and the sixth-order extrapolation method derived from Strang splitting method show good performance in the experiments with the periodic and Dirichlet boundary conditions respectively.
Efficient numerical algorithms for the  perturbed mechanical systems are proved very useful for the model of turbulence  of fluids.  Other initial boundary value problems including Burgers' type non linearities could be integrated by higher order splitting procedures.

\end{document}